\documentclass{compositio}
\newif\ifpdf
\ifx\pdfoutput\undefined
	\pdffalse
\else
	\pdfoutput=1
	\pdftrue
\fi
\usepackage{amsmath}
\usepackage{amsfonts}
\usepackage{enumerate}
\usepackage{graphpap}
\usepackage{mathrsfs}
\ifpdf
\DeclareGraphicsExtensions{.pdf}
\usepackage{hyperref}
\else
	\DeclareGraphicsExtensions{.eps}
	\usepackage{hyperref}
\fi

\newcommand\bc{{\mathbb C}}
\newcommand\bq{{\mathbb Q}}
\newcommand\bp{{\mathbb P}}

\newcommand\bz{{\mathbb Z}}
\newcommand\bs{{\mathbb S}}
\newcommand\bff{{\mathbb F}}
\newcommand\bbb{{\mathbb B}}
\newcommand\bdy{{\mathbf y}}
\newcommand\scl{{\mathscr L}}
\newcommand\scm{{\mathscr M}}
\newcommand\scn{{\mathscr N}}
\newcommand\scb{{\mathscr B}}
\newcommand\scc{{\mathscr C}}
\newcommand\scp{{\mathscr P}}
\newcommand\sch{{\mathscr H}}

\DeclareMathOperator\gl{GL}
\DeclareMathOperator\aut{Aut}
\DeclareMathOperator\pgl{PGL}

\newtheorem{thm}{Theorem}[section]

\newtheorem{prop}[thm]{Proposition}
\newtheorem{cor}[thm]{Corollary}
\newtheorem{lema}[thm]{Lemma}

\newtheorem{paso}{Step}
\theoremstyle{remark}
\newtheorem{obs}[thm]{Remark}

\theoremstyle{definition}
\newtheorem{dfn}[thm]{Definition}
\newtheorem{ntc}[thm]{Notation}

\newtheorem{ejm}[thm]{Example}

\newtheorem{cnt}[thm]{Construction}

\begin{document}
\begin{abstract}
We prove the existence of complexified 
real arrangements with the same combinatorics but different embeddings in $\bp^2$.
Such pair of arrangements has an additional property: they admit conjugated
equations on the ring of polynomials over $\bq(\sqrt{5})$.
\end{abstract}
\title{Topology and combinatorics of real line arrangements}
\shortauthors{E. Artal, J. Carmona, J.I. Cogolludo and M. Marco}
\author[E. Artal]
{Enrique ARTAL BARTOLO}
\email{artal@unizar.es}
\address{Departamento de Matem\'aticas\\
Campus Plaza de San Francisco s/n\\
E-50009 Zaragoza SPAIN}

\author[J. Carmona]{Jorge CARMONA RUBER}
\email{jcarmona@sip.ucm.es}
\address
{Departamento de Sistemas inform\'aticos y programaci\'on,
Universidad Complutense,\newline
Ciudad Universitaria s/n
E-28040 Madrid SPAIN}

\author[J.I. Cogolludo]{Jos\'e Ignacio COGOLLUDO AGUST\'IN}
\email{jicogo@unizar.es}
\address{Departamento de Matem\'aticas\\
Campus Plaza de San Francisco s/n\\
E-50009 Zaragoza SPAIN}

\author[M. Marco]{Miguel MARCO BUZUN\'ARIZ}
\email{mmarco@unizar.es}
\address{Departamento de Matem\'aticas\\
Campus Plaza de San Francisco s/n\\
E-50009 Zaragoza SPAIN}
\dedication{Accepted in Compositio Math.}
\date\today

\thanks{First, third and fourth authors are partially supported by
BFM2001-1488-C02-02; second author is partially
supported by BFM2001-1488-C02-01}
\maketitle

\section{Introduction}

In a famous preprint \cite{ry:98}, G.~Rybnikov proves the existence of two line arrangements 
$\scl_1,\scl_2$ in $\bp^2:=\bc \bp^2$ which have the same combinatorics but whose pairs,
$(\bp^2,\bigcup\scl_1)$ and $(\bp^2,\bigcup\scl_2)$, are not homeomorphic. Unfortunately,
this work has not been published in its final form. This may be due to the difficulty of 
verifying some of its statements. In another work~(\cite{accm:03a}), some parts of 
Rybnikov's paper are presented and explained, including a detailed verification of the 
quoted result. Rybnikov also gives an outline for proving a stronger result: the 
fundamental groups of $\bp^2\setminus\bigcup\scl_1$ and $\bp^2\setminus\bigcup\scl_2$ are 
not isomorphic. Providing a complete proof of this result would be most interesting.

In the present work we are interested in the topology of complexified real arrangements.
Complexified real arrangements form a very important class of arrangements where more 
topological properties are known. In this work we prove the existence of complexified 
real arrangements with the same combinatorics but different embeddings in $\bp^2$.
It is worth mentioning that Rybnikov's arrangements do not admit real equations. 
The counterexamples obtained in this work have an additional property: they admit conjugated
equations on the ring of polynomials over a number field, namely $\bq(\sqrt{5})$. The goal of 
this paper is to show the existence of such a pair of arrangements. In the following outline
of the paper, we will stress the similarities and major diferences with Ribnikov's approach:
\begin{itemize}

\item 
The first step towards finding such pairs of arrangements involves finding combinatorics 
whose moduli space is not connected. As a first approach, we will find \emph{ordered} 
combinatorics whose moduli space is not connected (an easier task). We proceed by 
considering an ordered line arrangement with equations in a finite field and then 
studying its complex realizations. In Rybnikov's paper, the MacLane arrangement in 
$\bff_3 \bp^2$ was considered.

\item 
In a second step we exhibit two ordered line arrangements $\scm_1,\scm_2$ which have the 
same ordered combinatorics but whose pairs, $(\bp^2,\bigcup\scm_1)$ and $(\bp^2,\bigcup\scm_2)$, 
are not homeomorphic by an order-preserving map. In order to prove this result we need to obtain 
some non-generic braid monodromies. In Rybnikov's paper, the author constructs an invariant of the 
isomorphism class of fundamental groups of the complements of lines whose ordered combinatorics 
induce the \emph{identity} on the first homology group. In~\cite{accm:03a} we describe an extension 
of Rybnikov's invariant which has been implemented using Maple. Using it, we are able to 
double check Rybnikov's example. Unfortunately this invariant does not give any valuable
information in our example.

\item 
Finally, using the aforementioned arrangements $\scm_1, \scm_2$ we construct a pair of 
arrangements $\scl_1, \scl_2$ so that the pairs $(\bp^2,\bigcup\scl_1)$ and $(\bp^2,\bigcup\scl_2)$ 
are not homeomorphic. This idea is part of Rybnikov's strategy: the fact that $\scm_1,\scm_2$ are 
\emph{real} arrangements allows us to use simpler arguments.

\end{itemize} 

Note that in our case the fundamental groups of $\bp^2\setminus\bigcup\scl_1$ and 
$\bp^2\setminus\bigcup\scl_2$ have isomorphic profinite completions, but we do not 
know whether or not they are isomorphic.

\section{An arrangement over a finite field}
\label{finito}
In what follows, a brief description of combinatorics and ordered combinatorics of a line
arrangement will be given.

\begin{dfn}
A \emph{line combinatorics} (\emph{projective configuration} in \cite{ry:98})
is a triple $(\scl,\scp,\in)$ where $\scl$ and $\scp$ are finite subsets and
$\in$ is a relation between $\scp$ and $\scl$ satisfying:

\begin{itemize}
\item $\forall \ell,\ell'\in\scl$, $\ell\neq \ell'$, $\exists! \ 
p \in \scp$ such that $p\in \ell$ and $p\in \ell'$ (in notation $p=\ell \cap \ell'$). 

\item  $\forall p\in\scp$ $\exists \ \ell,\ell'\in\scl$, $\ell\neq \ell'$,
such that $p=\ell\cap \ell'$.
\end{itemize}
Analogously we define an \emph{ordered line combinatorics}, where the elements of
$\scl=\{\ell_0,\ell_1,\dots,\ell_n\}$ are ordered. 
If the context is not ambiguous, we will omit both $\scp$ and $\in$.
The cardinal $\# \{\ell \mid p\in \ell\}$ is denoted as $\nu (p)$ and referred 
to as the \emph{multiplicity} of $p$ (as usual, points of multiplicity 2 are also called
\emph{double points} and points of multiplicity 3 \emph{triple points}). The second condition in 
a line combinatorics $(\scl,\scp,\in)$ implies that $\nu (p)\geq 2 \ \ \forall p\in \scp$.
\end{dfn}

A simple way to obtain line combinatorics is via arrangements in $\bff_q\bp^2$ for some 
finite field $\bff_q$. In such cases it is more practical to consider arrangements of points 
and then dualize.

For example, the starting point in \cite{ry:98} is the MacLane arrangement. This arrangement
can be defined as follows. Consider $\bff_3\bp^2$ as the union of $\bff_3^2$ and the line 
at infinity. Consider also the points of $\bff_3^2\setminus\{0\}$: the MacLane combinatorics 
is the abstract line combinatorics corresponding to the dual of this $8$ point arrangement.

\begin{dfn} 
Let $\scl$ be a line combinatorics. An \emph{automorphism} of $\scl$ is a permutation of 
$\scl$ preserving the incidence relations. The set of such automorphisms is the 
\emph{automorphism group} of $\scl$.
\end{dfn} 

\begin{ejm}It is easily seen that the automorphism group of MacLane combinatorics is naturally 
isomorphic to $\gl(2;\bff_3)$.
\end{ejm}

The starting point of this work is a similar idea. Let us consider $\bff_4\bp^2$ as the union 
of $\bff_4^2$ and the line at infinity $L_\infty$. Let us consider
$$
A:=
\begin{pmatrix}
0&1\\
1&\zeta
\end{pmatrix}
\in GL(2;\bff_4),
$$
where $\zeta$ is any of the elements of $\bff_4\setminus\bff_2$
(and hence the other one is $\bar\zeta=1+\zeta=\zeta^{-1}=\zeta^2$).
Note that $A^5=I_2$. The set
$$
\scb:=L_\infty\cup\left\{A^j
\left.
\begin{pmatrix}
1\\
0
\end{pmatrix}
\right\rvert 0\leq j\leq 4\right\}=
L_\infty\cup\left\{
\begin{pmatrix}
1\\
0
\end{pmatrix},
\begin{pmatrix}
0\\
1
\end{pmatrix},
\begin{pmatrix}
1\\
\zeta
\end{pmatrix},
\begin{pmatrix}
\zeta\\
\zeta
\end{pmatrix},
\begin{pmatrix}
\zeta\\
1
\end{pmatrix}
\right\}
$$
contains ten points, the last five of which will be denoted by $P_1,\dots,P_5$.
Note that each point lies in a different line through the origin. Let us denote 
by $Q_j$ the point at infinity of the line passing through both the origin and 
$P_j$, $j=1,\dots,5$. The following is a list of all the lines of 
$\bff_4\bp^2$ that contain at least two points (after each equation, a list
of points in $\scb$ contained in the line is provided):
\begin{itemize}
\item The line $L_\infty$ containing $Q_1,\dots,Q_5$.
\item Lines through $Q_1$: $y=0$ ($P_1$), $y=z$ ($P_2,P_5$), $y=\zeta z$
($P_3,P_4$).
\item Lines through $Q_2$: $x=z$ ($P_1,P_3$), $x=0$ ($P_2$), $x=\zeta z$
($P_4,P_5$).
\item Lines through $Q_3$: $y=\zeta(x+z)$ ($P_1,P_5$), $y=\zeta x+z$ ($P_2,P_4$), $y=\zeta x$ ($P_3$).
\item Lines through $Q_4$: $y=x+z$ ($P_1,P_2$), $y=x+\bar\zeta z$ ($P_3,P_5$), $y=x$ ($P_4$).
\item Lines through $Q_5$: $y=\bar\zeta(x+z)$ ($P_1,P_4$), $y=\bar\zeta x+z$ ($P_2,P_3$), 
$y=\zeta x$ ($P_5$).
\end{itemize} 

Note that each pair $P_i,P_j$, $1\leq i<j\leq 5$ appears exactly once in the above list.
The line combinatorics $\scc$ we are interested in is dual to $\scb$. 
It contains the following ten lines:

\begin{itemize}
\item $L_1,\dots,L_5$ dual to $P_1,\dots,P_5$ resp. and
\item $M_1,\dots,M_5$ dual to $Q_1,\dots,Q_5$ resp.
\end{itemize}   

The intersection points of these lines are:

\begin{itemize}
\item $O:=M_1\cap\dots\cap M_5$, dual to the origin 
$\left( \begin{smallmatrix} 0\\ 0 \end{smallmatrix} \right)$.
\item $P_{1,\{2,5\}}$, $P_{1,\{3,4\}}$, $P_{2,\{1,3\}}$, $P_{2,\{4,5\}}$, $P_{3,\{1,5\}}$, $P_{3,\{2,4\}}$, 
$P_{4,\{1,2\}}$, $P_{4,\{3,5\}}$, $P_{5,\{1,4\}}$, $P_{5,\{2,3\}}$, where $P_{i,\{j,k\}}:=M_i\cap L_j\cap L_k$.
\item $R_i:=M_i\cap L_i$, $i=1,\dots,5$, dual to the line passing through the origin,
$P_i$ and $Q_i$.
\end{itemize}

Note that the Falk-Sturmfels arrangement (see~\cite{Cohen-Suciu-braidmonod}) is the result of
deleting one of the $M_i$ lines from $\scc$.

\begin{lema}
\label{perm} 
Let $\psi\in\aut(\scc)$. Then there exists a unique $\sigma\in\Sigma_5$ 
such that $\psi(L_i)=L_{i^\sigma}$ and $\psi(M_i)=M_{i^\sigma}$,
$i=1,\dots,r$. The mapping $\psi\mapsto\sigma$ defines a monomorphism
$\aut(\scc)\to\Sigma_5$ which will be identified with the inclusion.
\end{lema}

\begin{proof}
First of all note that any $\psi\in\aut(\scc)$ preserves the multiplicity of any point.
Since $\{O\}=\{P\in \scc \mid \nu (P)=5\}$ one concludes that $\psi$ fixes $O$ and
hence should interchange the five lines $M_1,...,M_5$ containing $O$.
Let us denote by $\sigma\in\Sigma_5$ the permutation satisfying $\psi(M_i)=M_{i^\sigma}$.
Note that each line $M_i$ determines $R_i$ (as the only point of multiplicity 2
contained in $M_i$) and hence it determines $L_i$ (as the other line containing
$R_i$). Therefore $\psi (R_i)=R_{i^\sigma}$ and $\psi (L_i)=L_{i^\sigma}$.
\end{proof}

\begin{obs} 
Because of this identification $\aut(\scc)$ can be considered to act on $\scc$ on the right.
\end{obs}

\begin{ejm} The action of the matrix $A$ induces an automorphism 
$\psi_A$ of $\scc$ given by the $5$-cycle $(1,2,3,4,5)$.
\end{ejm}

\begin{ejm} Let us consider the automorphism of $\bff_4^2$ defined by
$B \sigma$ where 
$$B:=
\begin{pmatrix}
1&\zeta\\
0&\zeta
\end{pmatrix}
\in GL(2;\bff_4),$$
and $\sigma$ is the Galois involution of $\bff_4|\bff_2$.
The projective completion of this mapping preserves $\scb$ and induces an 
automorphism $\bar \psi_B$ of $\scc$ determined by the permutation $(2,4,5,3)$.
\end{ejm}

\begin{obs} Regarding again $\aut(\scc)$ as a subgroup of $\Sigma_5$, it is easily seen that 
$\bar \psi_B^{-1} \psi_A \bar \psi_B=\psi_A^3$, that is, 
$\langle \psi_A \rangle$ is a normal subgroup in $G:=\langle \psi_A,\bar \psi_B \rangle$. 
Since $G/\langle \psi_A \rangle\cong\langle \bar \psi_B \rangle$, the subgroup $G$ 
is a semidirect product of $\langle \psi_A \rangle\cong \bz_5$ and 
$\langle \bar \psi_B \rangle\cong \bz_4$ and hence has order $20$.
\end{obs}

\begin{lema} Under the above considerations $\aut(\scc)=G$.
\end{lema}

\begin{proof}
It is enough to check that $\aut(\scc)\subset G$. For the sake of simplicity, we will denote 
an automorphism in $\aut(\scc)$ and its induced permutation in $\Sigma_5$ by the same symbol.

Let us consider $\psi$ an automorphism fixing at least two lines of type $M_i$. 
After conjugation by an element in $G$, one may assume that these lines are $M_1$, $M_2$. 
Hence $L_1$ and $L_2$ are also fixed (Lemma~\ref{perm}). Using 
$\psi(P_{1,\{2,5\}})=P_{1,\{2,5^\psi\}}$, $\psi(P_{4,\{1,2\}})=P_{4^\psi,\{1,2\}}$, 
and $\psi(P_{2,\{1,3\}})=P_{2,\{1,3^\psi\}}$, the lines $L_5$, $L_4$, and $L_3$ are also 
respectively fixed, hence $\psi$ is the identity. 

Let us consider $\psi$ an automorphism fixing only one line of type $M_i$. Again, after 
conjugation by an element in $G$ one may suppose such a line is $M_1$. Since $\bar \psi_B$ fixes $M_1$ and 
permutes the others, after multiplying $\psi$ by a certain power of $\bar \psi_B$ one may 
suppose that $M_2$ is also fixed; then by the previous step $\psi \in G$.

Finally, let us consider a general $\psi$. Since $\psi_A$ is transitive, after multiplying 
$\psi$ by a certain power of $\psi_A$ one may suppose that $M_1$ is fixed; then by the 
previous step $\psi \in G$.
\end{proof}

\begin{obs} 
Following the notation of the previous examples, the automorphisms $\psi_M$
defined by a linear transformation $M$ will be called \emph{positive}, whereas the 
automorphisms $\bar \psi_M$ defined by $M \sigma$ ($\sigma$ being the Galois transformation 
of $\bff_4| \bff_2$) will be called \emph{negative}. When the automorphisms are viewed
as permutations, the concepts of positive and negative agree with their signature.
\end{obs}

\section{Complex realization}
\begin{dfn} 
Let $\scc$ be a line combinatorics. An arrangement $\scl$ of $\bp^2=\bc\bp^2$ is a 
\emph{complex realization of $\scc$} if its combinatorics agrees with $\scc$. An 
\emph{ordered complex realization} of an ordered line combinatorics is defined accordingly.
\end{dfn}

\begin{ntc} The space of all complex realizations of a line combinatorics $\scc$ is denoted 
by $\Sigma(\scc)$. This is a quasiprojective subvariety of $\bp^{\frac{n(n+3)}{2}}$, where
$n:=\#\scc$. If $\scc$ is ordered, we denote by $\Sigma^{\text{ord}}(\scc)\subset(\check\bp^2)^n$ 
the space of all ordered complex realizations of $\scc$.
\end{ntc}

There is a natural action of $\pgl(3;\bc)$ on such spaces. This justifies the following definition.

\begin{dfn} The \emph{moduli space} of a combinatorics $\scc$
is the quotient $\scm(\scc):=\Sigma(\scc)/\pgl(3;\bc)$.
The \emph{ordered moduli space} $\scm^{\text{ord}}(\scc)$ of an 
ordered combinatorics $\scc$ is defined accordingly.
\end{dfn}

\begin{ejm} Let us consider the MacLane line combinatorics $\scn$.
It is well known that $\#\scm(\scn)=1$ and that $\#\scm^{\text{ord}}(\scn)=2$. 
One can find representatives having equations in the polynomial ring 
over the field of cubic roots of unity. Moreover, the MacLane line combinatorics has no real
realization.
\end{ejm}

\begin{obs} Let us consider an ordered line combinatorics $\scc$. It is easily seen that 
$\aut(\scc)$ acts on both $\scm(\scc)$ and $\scm^{\text{ord}}(\scc)$. 
\end{obs}

\begin{ejm} The action of $\aut(\scn)\cong\gl(2,\bff_3)$ on the moduli spaces 
is as follows: matrices of determinant $+1$ (resp. $-1$) fix (resp. interchange) 
the two elements of $\scm^{\text{ord}}(\scn)$.
Of course complex conjugation also acts on $\scm^{\text{ord}}(\scn)$ interchanging the 
two elements. From the topological point of view, after fixing two representatives 
$\scn^\pm$ of $\scn$ one has that:
\begin{itemize}
\item There exists a homeomorphism $(\bp^2,\bigcup\scn^+)\to(\bp^2,\bigcup\scn^-)$
preserving orientations on both $\bp^2$ and the lines. Such a homeomorphism does not 
respect the ordering.

\item There exists a homeomorphism $(\bp^2,\bigcup\scn^+)\to(\bp^2,\bigcup\scn^-)$ preserving 
orientations on $\bp^2$, but not on the lines. Such a homeomorphism respects the ordering.
\end{itemize}
\end{ejm}

In what follows, we study the space and moduli of complex realizations for the
combinatorics $\scc$ introduced in \S\ref{finito}. Let us fix an ordering in the combinatorics 
so as to choose coordinates. Let us suppose that a complex realization exists. 
Up to the action of $\pgl(3;\bc)$ one may assume that:
$$
M_1: z=0,\quad
M_2: x=0,\quad
M_3: x=z,\quad
M_4: x=\alpha z,\quad
M_5: x=\beta z
$$
where $\alpha,\beta\in\bc\setminus\{0,1\}$, $\alpha\neq\beta$ (note that $\alpha,\beta$ 
cannot be fixed a priori since they represent cross-ratios, which are projective invariants).
Taking the incidence relations into account and by means of the subgroup of projective 
transformations which fix $M_1$, $M_2$, and $M_3$, one can also assume that $L_5: y=0$ and 
$L_4: y=z$. Note that the only projective transformation which fixes $M_1$, $M_2$, $M_3$, $L_4$, 
and $L_5$ is the identity. Also note that:

\begin{itemize}
\item $L_1$ passes through $L_5\cap M_2=[0:0:1]$
and $L_4\cap M_3=[1:1:1]$, hence $L_1: x=y$.

\item $L_2$ passes through $L_4\cap M_5=[\beta:1:1]$,
$L_5\cap M_3=[1:0:1]$, and $L_1\cap M_4=[\alpha:\alpha:1]$, hence:
$$
L_2: \alpha x-(\alpha-1) y-\alpha z=0,\quad
\beta=\frac{2\alpha-1}{\alpha}.
$$
\item $L_3$ passes through
$L_1\cap M_5=[2\alpha-1:2\alpha-1:\alpha]$,
$L_4\cap M_2=[0:1:1]$, $L_5\cap M_4=[\alpha:0:1]$,
and $L_2\cap M_1=[\alpha-1:\alpha:0]$. In principle one should obtain four
equations out of these conditions, but only three of them are independent, yielding:
$$
L_3:x+\alpha y-\alpha z=0,\quad
\alpha^2+\alpha-1=0.
$$
\end{itemize}

The family $\scm^{\text{ord}}(\scc)$ depends on one parameter $\alpha$ satisfying 
the above equation. The other root of this equation is $\gamma:=-1-\alpha$, which 
also verifies $\gamma^2+\gamma-1=0$, will be taken as the parameter.
Summarizing, one has the following:

\begin{prop}
\label{ecuaciones} 
The space $\scm^{\text{ord}}(\scc)$ has two elements. Moreover, representatives can be 
chosen to have the following equations:
$$
M_1: z=0,\quad
M_2: x=0,\quad
M_3: x=z,\quad
M_4: x=-(\gamma+1)z,\quad
M_5: x=(\gamma+2)z,
$$
$$
L_1: y=x,\quad
L_2: y=\gamma(x-z),\quad
L_3: y=\gamma x+z,\quad
L_4: y=z,\quad
L_5: y=0,
$$
where $\gamma^2+\gamma-1=0$.
\end{prop}

\begin{obs} 
The previous proposition implies that there are representatives of the two ordered complex 
realizations of $\scc$ having real equations conjugated in the ring of polynomials over the
number field $\bq(\sqrt{5})$. It is easily seen that the positive (resp. negative) automorphisms 
of $\scc$ preserve (resp. interchange) the elements in $\scm^{\text{ord}}(\scc)$. 
In particular, $\#\scm(\scc)=1$.
\end{obs}

\section{Topology of the realizations}

The first main result of Rybnikov's work~\cite{ry:98} reads as follows: let us 
denote by $\scn^\pm$ two representatives of the elements of $\scm^{\text{ord}}(\scn)$.
Then there is no isomorphism of the fundamental groups of $\bp^2\setminus\bigcup\scn^+$ and
$\bp^2\setminus\bigcup\scn^-$ inducing the identity in the homology groups. In particular, 
there is is no orientation-preserving homeomorphism between $(\bp^2,\bigcup\scn^+)$ and 
$(\bp^2,\bigcup\scn^-)$.

\begin{figure}[ht]
\centering
\includegraphics[height=5cm,width=5cm]{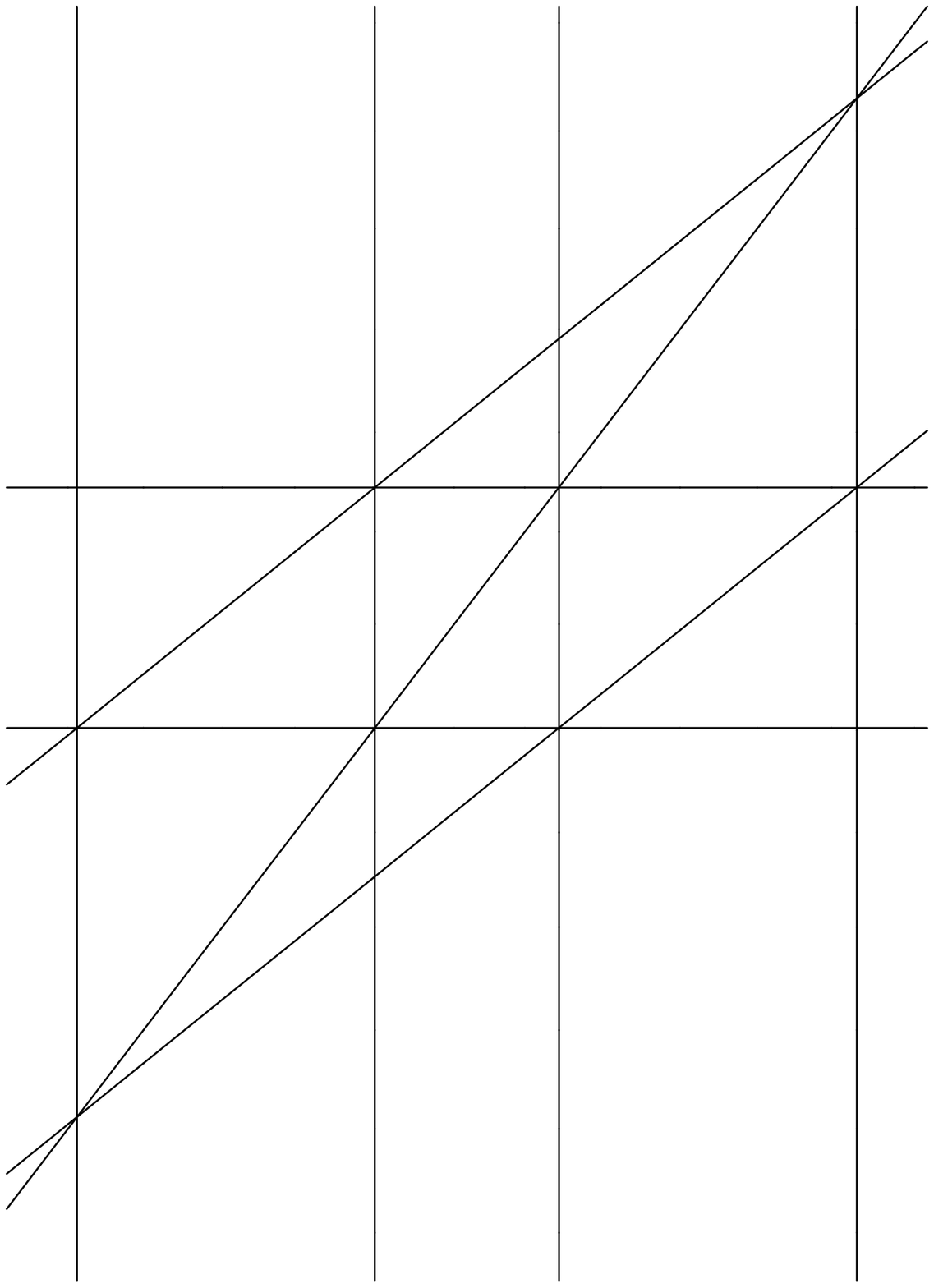} 
\caption{Arrangement $\scc^+$}
\label{recta1}
\begin{picture}(0,0)(0,3)
\put(-90,130){$L_4^+$}
\put(75,105){$L_5^+$}
\put(-88,85){$L_3^+$}
\put(37,140){$L_1^+$}
\put(75,150){$L_2^+$}
\put(-85,180){$M_4^+$}
\put(-34,180){$M_2^+$}
\put(17,60){$M_3^+$}
\put(63,60){$M_5^+$}
\end{picture}
\end{figure}

\begin{figure}[ht]
\centering
\includegraphics[height=5cm,width=5cm]{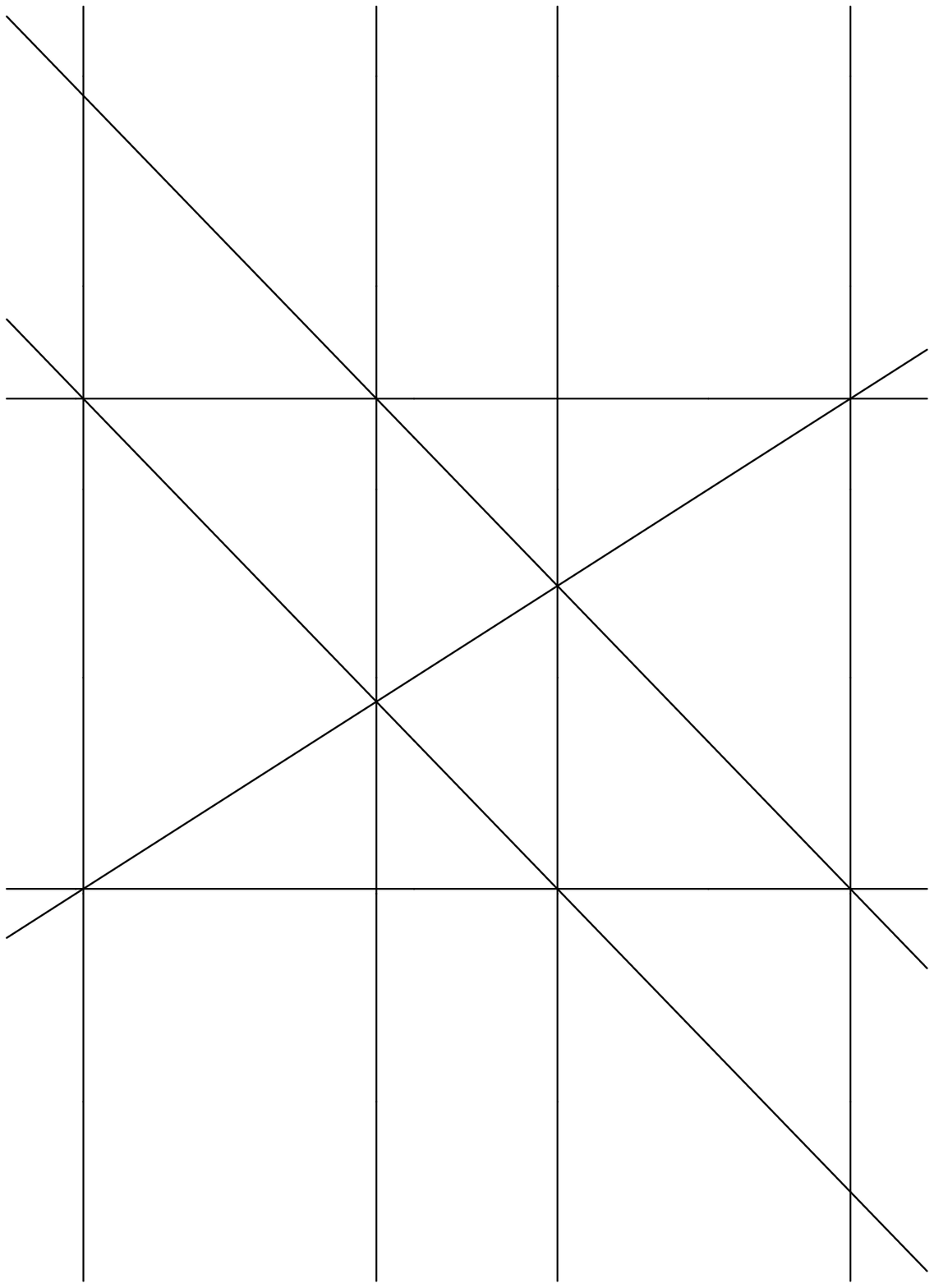} 
\caption{Arrangement $\scc^-$}
\label{recta2}
\begin{picture}(0,0)(0,3)
\put(30,149){$L_4^-$}
\put(-40,75){$L_5^-$}
\put(-83,158){$L_3^-$}
\put(35,117){$L_1^-$}
\put(75,70){$L_2^-$}
\put(15,180){$M_4^-$}
\put(-80,60){$M_2^-$}
\put(62,180){$M_3^-$}
\put(-35,55){$M_5^-$}
\end{picture}
\end{figure}

Analogously, one can consider the two ordered arrangements $\scc^\pm$ representing 
the elements of $\scm^{\text{ord}}(\scc)$. They are determined by the equations
given in Proposition~\ref{ecuaciones} and by the choice of 
$\gamma^\pm:=\frac{-1\pm\sqrt{5}}{2}$.
Let us denote by $L_i^\pm$ and $M_i^\pm$, $i=1,\dots,5$ the lines as in section \S 2.
Affine pictures are shown in Figures \ref{recta1} and \ref{recta2}, where $M_1^\pm$ 
are the corresponding lines at infinity.

As pointed out in the introduction, we have completely checked Rybnikov's computations 
concerning MacLane arrangements~\cite{accm:03a}. Rybnikov's approach does not work with $\scc$ 
and we do not know whether or not there exists an isomorphism of the fundamental groups of
$\bp^2\setminus\bigcup\scc^+$ and $\bp^2\setminus\bigcup\scc^-$ inducing the identity in 
the homology groups. Nevertheless, we are able to prove that there is no orientation-preserving 
homeomorphism $(\bp^2,\bigcup\scc^+)\to(\bp^2,\bigcup\scc^-)$. The required invariant is a 
non-generic braid monodromy. Let us introduce some concepts in order to adapt braid monodromy 
to line arrangements.

\begin{dfn} An \emph{affine horizontal arrangement} $\scl$
is an affine line arrangement in $\bc^2$ with no vertical lines. The ramification 
set $A_\scl\subset\bc$ is the set of $x$-coordinates of the multiple points of $\scl$.
\end{dfn}

\begin{dfn} An \emph{ordered affine horizontal arrangement} $\scl$
is an affine horizontal arrangement with an ordering on both $\scl$ and $A_\scl$. 
By means of the usual embedding $\bc^2\subset\bp^2$, the \emph{ordered fibered projective arrangement} $\scl^\varphi$ associated with $\scl$ is obtained by considering:
\begin{itemize}
\item the projectivization of the lines in $\scl$;
\item the line at infinity $L_\infty$;
\item the projectivization of the fibers of $A_\scl$ (the vertical lines, pre-images of $A_\scl$).
\end{itemize}
\end{dfn}

\begin{ejm} 
The non-vertical lines in Figures~\ref{recta1}~and~\ref{recta2} give two
affine horizontal arrangements $\scl^\pm$. With a suitable order,
their ordered fibered projective arrangements are $\scc^\pm$.
\end{ejm}

\begin{obs} Given an ordered affine horizontal arrangement $\scl$,
the mapping $\bp^2\setminus\scl^\varphi\to\bc\setminus A_\scl$, induced by the projection 
$(x,y)\mapsto x$ is a locally trivial fibration.
\end{obs}

\begin{dfn} Let $Z$ be a connected projective manifold
and let $H$ be a hypersurface of $Z$. Let $*\in Z\setminus H$
and let $K$ be an irreducible component of $H$. A homotopy
class $\gamma\in\pi_1(Z\setminus H;*)$ is called a
\emph{meridian about $K$ with respect to $H$}
if $\gamma$ has a representative $\delta$ satisfying the
following properties:
\begin{enumerate}[(a)]
\item there is a smooth complex analytic disk
$\Delta \subset Z$ transverse to $H$ such that
$\Delta\cap H=\{*'\} \subset K$
(transversality implies that $*'$ is a smooth point of $H$).

\item there is a path $\alpha$ in $Z\setminus H$ from $*$
to $*'' \in \partial\Delta$.

\item $\delta=\alpha\cdot\beta\cdot\alpha^{-1}$, where
$\beta$ is the closed path obtained by traveling from $*''$
along $\partial\Delta$ in the positive direction.
\end{enumerate}
\end{dfn}

\begin{dfn} Let $A:=\{a_1,\dots,a_r\}\subset\bc$ and let $\Delta\subset\bc$ be a disk 
containing $A$ in its interior. Fix $*\in\partial\Delta$. A \emph{geometric basis} of the 
free group $\pi_1(\bc\setminus A;*)$ is a basis $\gamma_1,\dots,\gamma_r$ where
each $\gamma_i$ is a meridian of a point in $A$ and $\gamma_r\cdot\ldots\cdot\gamma_1=\partial\Delta$
positively oriented.
\end{dfn}

\begin{obs} In order to compute braid monodromies
we use the \emph{reversed} lexicographic order on $\bc$, i.e., if
$z,w\in\bc$, $z\neq w$, then $z\succ w$ if
$\Re z>\Re w$ or $\Re z=\Re w$ and $\Im z>\Im w$.
\end{obs}

\begin{cnt} Let us construct a special {geometric basis} of $\pi_1(\bc\setminus A;*)$, where 
 $*$ is a \emph{big enough} real number and $A:=\{a_1,\dots,a_r\}\subset\bc$. Let us order the elements 
of $A$ such that $a_1\succ \dots\succ a_r$. Consider the polygonal path $\Gamma$ with vertices 
$*,a_1,\dots,a_r$. Let us choose a small enough $\varepsilon>0$ so that the disks $\Delta_i$ 
of radius $\varepsilon$ centered at $a_1,\dots,a_r$ are pairwise disjoint. The geometric 
basis $\gamma_1,\dots,\gamma_r$ is defined as follows. Fix $i\in\{1,\dots,r\}$; let us consider 
$\alpha_i$ a path starting at $*$ and running along $\Gamma$ outside the interior of the disks 
$\Delta_1,\dots,\Delta_{i-1}$ and counterclockwise along the arcs
$\partial\Delta_1,\dots,\partial\Delta_{i-1}$. This path should stop at the first intersection 
$*_i$ of $\Gamma$ with $\partial\Delta_i$. Let us denote by $\delta_i$ the loop based at 
$*_i$ and going counterclockwise along $\partial\Delta_i$. Then $\gamma_i:=\alpha_i\cdot\delta_i\cdot\alpha_i^-$. 
This basis will be called a \emph{lexicographic basis}.
\end{cnt}

\begin{ejm} Let us denote by $\alpha_5^+$, $\alpha_3^+$, $\alpha_2^+$, and $\alpha_4^+$ a 
lexicographic geometric basis of $\pi_1(\bc\setminus A_{\scl^+};R)$, $R\gg0$; the indices 
of $\alpha_i^+$ refer to the lines they are meridians of. Analogously, $\alpha_3^-$, 
$\alpha_4^-$, $\alpha_5^-$, and $\alpha_2^-$ will denote a lexicographic geometric basis 
of $\pi_1(\bc\setminus A_{\scl^-};R)$. 
Let us recall that $A_{\scl^\pm}=\{\gamma^\pm +2,1,0,-(\gamma^\pm +1)\}$.
\end{ejm}

Let $\bbb_n$ be the abstract braid group on $n=\# \scc$ 
strings given by the standard presentation:
$$
\left\langle \sigma_1,\dots,\sigma_{n-1} :\
[\sigma_i,\sigma_j]=1,\ |i-j|\geq 2,
\sigma_i\sigma_{i+1}\sigma_i=\sigma_{i+1}\sigma_i\sigma_{i+1},\
i=1,\dots,n-2
\right\rangle.
$$
The pure braid group $\bp_n$ is the kernel of the natural epimorphism of
$\bbb_n$ onto the symmetric group $\Sigma_n$ (obtained from the above presentation
by adding the relations $\sigma_i^2=1$).

\smallbreak
Let $\bdy\subset\bc$ be a subset of exactly $n$ elements. Let us recall that the braid group 
$\bbb_{\bdy}$ is the set of homotopy classes, relative to $\{0,1\}$, of sets of paths
$\{\gamma_1,\dots,\gamma_n\}$, $\gamma_j:[0,1]\to\bc$,
starting and ending at $\bdy$ and such that
$\forall t\in[0,1]$, $\{\gamma_1(t),\dots,\gamma_n(t)\}$
is a set of $n$ distinct points. The elements of
$B_{\bdy}$ are called braids based at $\bdy$ and
are represented as a set of non-intersecting paths in
$\bc\times[0,1]$, as usual. The pure braids correspond to sets of loops and they form the 
subgroup $\bp_\bdy$.
Products are also defined in the standard way.
Analogously, if $\bdy_1,\bdy_2\subset\bc$ are subsets of exactly $n$ elements then 
$\bbb_{\bdy_2,\bdy_1}$ is the set of braids starting at $\bdy_1$ and ending at $\bdy_2$. 
The abstract braid group is canonically identified with $\bbb_{\bdy_0}$, where 
$\bdy_0:=\{-1,\dots,-n\}$. The groups $\bp_{\bdy_0}$ and $\bp_n$ are identified as well.

\begin{ejm} 
There is a special element of $\bbb_{\bdy_2,\bdy_1}$
which is called the $(\bdy_2,\bdy_1)$-lexicographic braid. It is easily defined as follows.
Let us order $\bdy_1$ and $\bdy_2$ according to the lexicographic order. The braid is 
defined by taking the linear segments joining the corresponding elements in $\bdy_1$ and 
$\bdy_2$. The main point of this construction is that given $\bdy_1$ and $\bdy_2$ one can identify
$\bbb_{\bdy_2,\bdy_1}$ with $\bbb_n$ by means of the lexicographic braids associated with
$(\bdy_0,\bdy_1)$ and $(\bdy_2,\bdy_0)$.
\end{ejm}

The braid monodromy of an ordered affine horizontal arrangement $\scl$ can be defined as follows.
Let us choose $*\in\bc\setminus A_\scl$ and let $\bdy^*\subset\bc$ the set of $y$-coordinates 
of the points in $\scl$ with $x$-coordinate $*$. Moving around a loop in $\bc\setminus A_\scl$
based at $*$ one can obtain a braid in $\bbb_{\bdy^*}$ and conjugating by the  $(\bdy_0,\bdy^*)$-lexicographic braid one obtains a braid in $\bp_n$.

\begin{dfn} 
The \emph{braid monodromy} of an ordered affine horizontal arrangement $\scl$
is the group homomorphism $\rho:\pi_1(\bc\setminus A_\scl;*)\to\bp_n$ defined above.
\end{dfn}

\begin{obs} 
The braid monodromy is well defined up to conjugation in $\bp_n$.
\end{obs}

\begin{prop} 
The braid monodromies $\rho^\pm$ of the arrangements $\scl^\pm$ are given by:
\begin{equation*}
\begin{split}
\rho^+(\alpha_5^+)=\rho^-(\alpha_3^-)=&\sigma_1^2\sigma_3^2\\
\rho^+(\alpha_3^+)=\rho^-(\alpha_4^-)=&(\sigma_1\sigma_3)* (\sigma_2^2\sigma_4^2)\\
\rho^+(\alpha_2^+)=\rho^-(\alpha_5^-)=&(\sigma_1\sigma_3\sigma_2 \sigma_4)* (\sigma_1^2\sigma_3^2)\\
\rho^+(\alpha_4^+)=\rho^-(\alpha_2^-)=&(\sigma_1\sigma_3\sigma_2 \sigma_4\sigma_1\sigma_3)* 
(\sigma_2^2\sigma_4^2)
\end{split}
\end{equation*}
where $a*b:=aba^{-1}$ and the strings in the  positive 
(resp. negative) case correspond to the lines $(L_1^+,L_3^+,L_2^+,L_4^+,L_5^+)$, resp. 
$(L_1^-,L_4^-,L_5^-,L_2^-,L_3^-)$.
\end{prop}

\begin{proof} 
It is straightforward from Figures~\ref{recta1}~and~\ref{recta2} and the choice of 
lexicographic braids.
\end{proof}

\begin{obs}
\label{cambio} 
Replacing $\rho^-(\alpha_j^-)$ by $\tau^{-1}\rho^-(\alpha_j^-)\tau$, where
$\tau:=\sigma_3\sigma_4\sigma_2\sigma_3\sigma_2$,
the strings of the braids correspond to the \emph{same} sets of ordered lines. 
After this change, we will keep the original notation $\rho^-$.
\end{obs}

Our goal is to adapt~\cite[Theorem 1]{acc:01a} to the case of ordered horizontal 
line arrangements. The main difference is that the orbit of an ordered arrangement is
given by the action of the pure braid group instead of the whole braid group as is
required in the unordered case. Let us fix an ordered horizontal line arrangement $\scl$. Let 
us consider its braid monodromy $\rho:\pi_1(\bc\setminus A_\scl)\to\bp_n$. Fixing a 
geometric basis $\gamma_1,\dots,\gamma_r$ of $\pi_1(\bc\setminus A_\scl)$, $r=\# A_\scl$, 
the $r$-tuple $(\rho(\gamma_1),\dots,\rho(\gamma_r))$ is said to determine, or 
\emph{represent}, the braid monodromy $\rho$. Recall that $\rho$ was defined up 
to conjugation in $\bp_n$; this comes from the two choices that have been made: the 
base point $*$ and the geometric basis.

\begin{dfn} 
Two elements of $\bp_n^r$ are said to be \emph{equivalent}
if they can be related by simultaneous conjugation in $\bp_n$ and
\emph{Hurwitz moves}. Recall that the $j^{\text{th}}$ Hurwitz move
of $(\tau_1,\dots,\tau_r)\in\bp_n^r$ produces
$$
(\tau_1,\dots,\tau_{j-1},\tau_{j+1},\tau_{j+1}*\tau_j,\tau_{j+2},\dots,\tau_r).
$$
\end{dfn}

\begin{obs} 
Given $\tau:=(\tau_1,\dots,\tau_r)\in\bp_n^r$ it is important to consider its 
\emph{pseudo-Coxeter element} $\chi (\tau):=\tau_r\cdot\ldots\cdot\tau_1\in \bp_n$. 
Note that $\chi(\tau)$ remains invariant by Hurwitz moves.
\end{obs}

\begin{lema} 
Under the above notation, let $(\rho(\gamma_1),\dots,\rho(\gamma_r))$ represent the
braid monodromy of $\scl$. Then $\tau:=(\tau_1,\dots,\tau_r)\in\bp_n^r$
also represents the braid monodromy of $\scl$ if and only if $\tau$ is equivalent to 
$(\rho(\gamma_1),\dots,\rho(\gamma_r))$.
\end{lema}

\begin{proof}
The ideas involved in the proof of this lemma are classical, see \cite{acc:01a}.
\end{proof}

\begin{thm}
\label{top}
Let $\scl_1,\scl_2\subset\bp^2$ be two ordered affine horizontal line arrangements such 
that their associated fibered projective arrangements have the same combinatorics.
Let also $F:\bp^2\to\bp^2$ be an orientation-preserving homeomorphism such that:
\begin{enumerate}[{\rm(i)}]
\item $F(P)=P$, $P$ being the projection point,

\item $F(\scl_1)=\scl_2$ respecting ordering and orientations, and

\item $F$ respects the vertical lines preserving both ordering and orientations.
\end{enumerate}
Then both ordered affine horizontal line arrangements have equivalent braid monodromies.
\end{thm}

\begin{proof}
The proof essentially follows the one in~\cite[Theorem 1]{acc:01a}.
\end{proof}

\begin{thm}
\label{orden} 
There is no order-preserving homeomorphism between $(\bp^2,\bigcup\scc^+)$ and $(\bp^2,\bigcup\scc^-)$.
\end{thm}

\begin{proof}
We break the proof up into several steps.

\begin{paso}
\label{masmas} 
There is no homeomorphism between $(\bp^2,\bigcup\scc^+)$ and $(\bp^2,\bigcup\scc^-)$ 
preserving the ordering and the orientation on both the plane and the lines.
\end{paso}

\begin{proof}
This result follows from Theorem~\ref{top}. Recall that $\rho^-$ has been replaced by 
a conjugate (Remark~\ref{cambio}). Let us suppose that $\rho^+$ and $\rho^-$ are equivalent.

Let $K^\pm$ denote the image of $\rho^\pm$ in the pure braid group. Let $c^\pm:=\chi(\rho^\pm(\alpha^\pm))$ be the pseudo-Coxeter elements of
$$
\rho^+(\alpha^+):=(\rho^+(\alpha^+_5),\rho^+(\alpha^+_3),\rho^+(\alpha^+_2),\rho^+(\alpha^+_4))
\text{ and }
\rho^-(\alpha^-):=(\rho^-(\alpha^-_3),\rho^-(\alpha^-_4),\rho^-(\alpha^-_5),\rho^-(\alpha^-_2));
$$ 
since $c^\pm$ is invariant by Hurwitz moves, one should have that $c^+$ and $c^-$ are 
conjugate in the pure braid group $\bp_5$. Replacing $\rho^-$ by the action of an appropriate 
inner automorphism, one may assume that $c:=c^+=c^-$. Let $C_c$ be the centralizer of $c$ in
$\bp_5$. Since $\rho^+(\alpha^+)$ and $\rho^-(\alpha^-)$ are equivalent, $K^+$ and $K^-$ 
should be conjugate by an element $g\in \bp_5$ such that 
$g \chi(\rho^+(\alpha^+)) g^{-1}=\chi(\rho^-(\alpha^-))$, that is, $gc=cg$, hence $g\in C_c$
(note that Hurwitz moves preserve the groups $K^\pm$ and the pseudo-Coxeter elements).

Let us consider the Burau representation of $\bbb_5$ into $\gl(5;\bz[t^{\pm 1}])$. 
Replacing $t$ by $2\mod 5$ one obtains a representation $\beta:\bbb_5\to\gl(5;\bz/5\bz)$. 
The subgroup $\beta(\bp_5)$ has order $58032\cdot 10^6$. Let us denote 
$\tilde C_c:=\beta(C_c)$, and $\tilde K^\pm:=\beta(K^\pm)$. Checking that 
$\tilde K^+$ and $\tilde K^-$ are not conjugate in $\gl(5;\bz/5\bz)$ by an element in 
$\tilde C_c$ is a matter of comparing finite tables; the groups $\tilde K^\pm$ have order 
$30000$ while $\#\tilde C_c=115200$. This has been verified using a GAP4~\cite{GAP4} program.
Its source can be found in the appendix or downloaded from a public 
site\footnote{ \href{http://riemann.unizar.es/geotop/pub/gap-programs/lines/monodromy\_groups}{\texttt{http://riemann.unizar.es/geotop/pub/gap-programs/lines/monodromy\_groups}}}.
\end{proof}

\begin{paso} 
There is no homeomorphism between $(\bp^2,\bigcup\scc^+)$ and $(\bp^2,\bigcup\scc^-)$ 
preserving the ordering, the orientations on the plane, and reversing the orientations on the lines.
\end{paso}

\begin{proof}
Composing such a homeomorphism with the complex conjugation (which preserves the orientation
of the plane, but reverses the orientations on the lines) one obtains a homeomorphism 
that is ruled out by Step~\ref{masmas}.
\end{proof}

For the remaining cases the following lemma is required.

\begin{lema}
\label{linking} 
Fix an orientation on $\bs^3$ and consider $\varphi:\bs^3\to\bs^3$ a homeomorphism, 
$K_1,K_2,L_1,L_2\subset\bs^3$ oriented knots such that $K_1\cap K_2=\emptyset$ and
$\varphi(K_i)=L_i$ (disregarding orientations), $i=1,2$.
Let us set $\vartheta$ and $\varepsilon_i=\pm 1$, $i=1,2$, such that:
$$
\vartheta:=
\begin{cases}
1&\text{if $\varphi$ preserves orientation}\\
-1&\text{otherwise,}
\end{cases}
\qquad
\quad
L_i=\varepsilon_i K_i.
$$
Then, if $\scl_{\bs^3}$ is the linking number, 
$\scl_{\bs^3}(L_1,L_2)=
\vartheta\varepsilon_1\varepsilon_2\scl_{\bs^3}(K_1,K_2)$.
\end{lema}

\begin{obs}
\label{linkingnumber}
Note also that, from the well-known relationship between linking and intersection numbers, 
if two lines intersect at one point and one considers the oriented knots they define 
on the oriented boundary of a small ball centered at the intersection point, then their 
linking number equals $1$.
\end{obs}

\begin{paso} 
There is no homeomorphism between $(\bp^2,\bigcup\scc^+)$ and $(\bp^2,\bigcup\scc^-)$ 
preserving the ordering, the orientation on the plane, and the orientation on at least one line.
\end{paso}

\begin{proof}
By Step~\ref{masmas} there should be at least one line whose orientation is reversed and by 
hypothesis, at least one line whose orientation is preserved. Both lines should of
course intersect at one point. Let us consider $K$ and $K'$ their links at the 
intersection point. According to the notation of Lemma~\ref{linking}, the hypotheses
imply that $\vartheta=1$ and $\varepsilon\varepsilon'=-1$, and hence $\scl_{\bs^3}(L,L')=-1$ 
where $L$ and $L'$ are the links of two lines around their intersection point. 
This contradicts Remark~\ref{linkingnumber}.
\end{proof}

\begin{paso} 
There is no homeomorphism between $(\bp^2,\bigcup\scc^+)$ and $(\bp^2,\bigcup\scc^-)$ 
preserving the ordering and reversing the orientation on the plane.
\end{paso}

\begin{proof}
Since $\# \scc^+>2$, there exist at least two lines whose orientations are simultaneously 
either preserved or reversed. Consider their links at the intersection point. According to 
the notation of Lemma~\ref{linking}, the hypotheses imply that $\vartheta=-1$ and 
$\varepsilon\varepsilon'=1$, which again contradicts Remark~\ref{linkingnumber}.
\end{proof}

This concludes the result.
\end{proof}

\section{A counterexample}

Since the aforementioned line combinatorics $\scn$ and $\scc$ have a great symmetry, they 
verify the equalities $\#\scm(\scn)=\#\scm(\scc)=1$ and hence will not provide pairs of arrangements 
having the same combinatorics but different topologies. In \cite{ry:98} the problem is 
solved by \emph{juxtaposition} of two arrangements with combinatorics of type $\scn$. 
In this work we will follow a simpler approach than can be applied in many different ways. 
For the sake of brevity, we describe only one example.

Let $\sch^\pm:=\scc^\pm\cup\{N^\pm\}$ be two new arrangements, where $N^\pm$ is the line 
joining the points $L_3^\pm\cap L_5^\pm\cap M_4^\pm=[1:0:-\gamma^\pm]$ and
$L_2^\pm\cap M_2^\pm=[0:1:-(\gamma^\pm+1)]$. In particular,
$N^\pm: \gamma^\pm x+(\gamma^\pm+1) y+z=0$. Note that these two arrangements
have the same combinatorics, say $\sch$, since their equations are conjugate in $\bq(\sqrt{5})$.
The new abstract line of $\sch$ will be denoted by $N$.

\begin{prop}
\label{combi-m} 
The automorphism group of $\sch$ is trivial.
\end{prop}

\begin{proof} 
Let $\psi\in\aut(\sch)$. Note that $\psi(N)\neq M_j$, $j=1,\dots,5$, since $N$ does not 
pass through the unique point of multiplicity $5$. Also, note that $N$ contains $3$ double 
points and each $L_j$ contains at most $2$ double points. Hence, $\psi(N)=N$,
$\psi(\scc)=\scc$ and therefore $\psi$ induces an automorphism of $\scc$. By simple 
combinatorial arguments, $\psi(M_4)=M_4$ and $\psi(M_2)=M_2$, therefore $\psi=1_\sch$. 
\end{proof}

\begin{cor} $\#\scm(\sch)=2$ and $\sch^\pm$ are representatives of the two elements.
\end{cor}

\begin{thm} There is no homeomorphism between $(\bp^2,\bigcup\sch^+)$ and $(\bp^2,\bigcup\sch^-)$.
\end{thm}

\begin{proof}
By Proposition~\ref{combi-m}, such a homeomorphism would induce the identity on $\sch$,
hence having to send $N^+$ to $N^-$. Therefore, it would define an order-preserving 
homeomorphism between $(\bp^2,\bigcup\scc^+)$ and $(\bp^2,\bigcup\scc^-)$. The result follows from 
Theorem~\ref{orden}.
\end{proof}

\begin{obs} 
An important feature of this pair of arrangements
is the fact that they are defined by conjugate equations in $\bq(\sqrt{5})$. This fact 
implies that they cannot be distinguished by algebraic methods. Note for example that 
the fundamental groups of $\bp^2\setminus\bigcup\sch^+$ and $\bp^2\setminus\bigcup\sch^-$ have 
isomorphic profinite completions.
\end{obs}

\section{Appendix}

In this appendix, we provide the source (with comments) of the short GAP4 program
mentioned in Step~\ref{masmas} of Theorem~\ref{orden}.
\vspace{3mm}
\hrule
\vspace{3mm}
\noindent\emph{\# $n$ is the number of strings.}\\
\texttt{n:=5;;}\\
\emph{\# The braid group in $n$ strings is constructed by generators and relations.}\\
\texttt{f:={\bf FreeGroup}(n-1);\\
x:={\bf GeneratorsOfGroup}(f);\\
rel:=[];\\
{\bf for} i in [1..n-3] {\bf do}\\
\hspace*{1cm} {\bf for} j in [i+2..n-1] {\bf do}\\
\hspace*{2cm}{\bf Add}(rel,{\bf Comm}(x[i],x[j]));\\
\hspace*{1cm}{\bf od};\\ 
{\bf od};\\
{\bf for} i in[1..n-2] {\bf do}\\
\hspace*{1cm}{\bf Add}(rel,x[i]*x[i+1]*x[i]/x[i+1]/x[i]/x[i+1]);\\
{\bf od};\\
g:=f/rel;\\
y:={\bf GeneratorsOfGroup}(g);\\}
\emph{\# The following is a system of generators of the pure braid group.}\\
\texttt{sub:=[y[1]\^{}2,y[2]\^{}2,y[1]*y[2]\^{}2*y[1],%
y[3]\^{}2,y[2]*y[3]\^{}2*y[2],\\
y[1]*y[2]*y[3]\^{}2*y[2]*y[1],
y[4]\^{}2,y[3]*y[4]\^{}2*y[3],\\
y[2]*y[3]*y[4]\^{}2*y[3]*y[2],
y[1]*y[2]*y[3]*y[4]\^{}2*y[3]*y[2]*y[1]];;\\
h:={\bf Subgroup}(g,sub);;\\}
\emph{\# The generators of $K^+$ will be defined, as well as the pseudo-Coxeter element 
and the braid used in Remark~{\rm\ref{cambio}} to obtain $K^-$ are given.}\\
\texttt{ sub1:=[y[1]\^{}2*y[3]\^{}2,(y[2]\^{}2*y[4]\^{}2)\^{}(y[1]\^{}-1*y[3]\^{}-1),\\
(y[1]\^{}2*y[3]\^{}2)\^{}(y[2]\^{}-1*y[4]\^{}-1*y[1]\^{}-1*y[3]\^{}-1),\\
(y[2]\^{}2*y[4]\^{}2)\!\^{}\!(y[1]\^{}-1*y[3]\^{}-1*y[2]\^{}-1*y[4]\^{}-1*y[1]\^{}-1*y[3]\^{}-1)\!]\!;\!;\\
del:=sub1[4]*sub1[3]*sub1[2]*sub1[1];;\\
cj:=y[3]*y[4]*y[2]*y[3]*y[2];\\
k1:={\bf Subgroup}(g,sub1);;}\\
\emph{\# The function} \texttt{burau}\emph{ is used to construct the composition of the Burau representation with 
the mapping $\bz[t^{\pm 1}]\to\bff_m$ such that $t$ is sent to the $n^{\text{th}}$ power 
of a generator of the invertible elements. 
The output of this function is a list of two elements containing 
the representation itself and the image of the pure braid group.}\\
\texttt{
burau:={\bf function}(n,m)\\
\hspace*{5mm}{\bf local} uno,cero,z,gg,hom,hh,s,R;\\
\hspace*{5mm}R:={\bf GF}(m);\\
\hspace*{5mm}s:={\bf Z}(m)\^{}n;;\\
\hspace*{5mm}uno:={\bf Z}(m)\^{}0;\\
\hspace*{5mm}cero:=0*{\bf Z}(m);\\ \hspace*{5mm}z:=[[[uno-s,s,cero,cero],[uno,cero,cero,cero],\\
\hspace*{5mm}[cero,cero,uno,cero],[cero,cero,cero,uno]],\\
\hspace*{5mm}[[uno,cero,cero,cero],[cero,uno-s,s,cero],\\
\hspace*{5mm}[cero,uno,cero,cero],[cero,cero,cero,uno]],\\      
\hspace*{5mm}[[uno,cero,cero,cero],[cero,uno,cero,cero],\\
\hspace*{5mm}[cero,cero,uno-s,s],[cero,cero,uno,cero]],\\
\hspace*{5mm}[[uno,cero,cero,-uno],[cero,uno,cero,-uno],\\
\hspace*{5mm}[cero,cero,uno,-uno],[cero,cero,cero,-s]]];\\
\hspace*{5mm}gg:={\bf Group}(z);\\
\hspace*{5mm}hom:={\bf GroupHomomorphismByImages}(g,gg,y,z);\\
\hspace*{5mm}hh:={\bf Image}(hom,h);\\
\hspace*{5mm}{\bf return} [hom,hh];\\
end;}\\
\emph{\# The representation $\beta$ is set in} \texttt{hom.}\\
\texttt{aux:=burau(1,5);\\
hom:=aux[1];\\
hh:=aux[2];;}\\
\emph{\# } \texttt{del1} \emph{and}\texttt{ del2} \emph{are the images of the pseudo-Coxeter elements.}\\
\texttt{del1:={\bf Image}(hom,del);\\
cj1:={\bf Image}(hom,cj);\\
del2:=del1\^{}cj1;}\\
\emph{\# The image} \texttt{nv} \emph{of a pure braid such that} \texttt{del1=del2\^{}nv} 
\emph{is obtained. We compute the group $\tilde K^+$;
the group $\tilde K^-$ is obtained after a conjugation
in order 
to have the equality of pseudo-Coxeter elements.}\\
\texttt{nv:={\bf RepresentativeAction}(hh,del2,del1);;\\
kk1:={\bf Image}(hom,k1);\\
kk2:=kk1\^{}(cj1*nv);;}\\
\emph{\# The centralizer of the pseudo-Coxeter element is computed and $\tilde K^\pm$ are 
verified to not be conjugated by an element in this centralizer.}\\
\texttt{hh1:=Centralizer(hh,del1);\\
IsConjugate(hh1,kk1,kk2);}
\vspace{3mm}
\hrule
\providecommand{\bysame}{\leavevmode\hbox to3em{\hrulefill}\thinspace}
\providecommand{\MRhref}[2]{%
  \href{http://www.ams.org/mathscinet-getitem?mr=#1}{#2}
}
\providecommand{\href}[2]{#2}
\providecommand\MR[1]{\relax\ifhmode\unskip\spacefactor3000 \space\fi
  \MRhref{#1}{#1}}

\end{document}